\documentclass[12pt,a4paper]{article}

\usepackage{graphicx}
\usepackage{amsmath}
\usepackage{amsfonts}
\usepackage{mathrsfs}
\usepackage{bbold}

\newtheorem{theorem}{Theorem}
\newtheorem{lemma}{Lemma}

\let\ds\displaystyle
\let\scr\mathscr
\newcommand\UU{{\mathbb U}}
\newcommand\CC{{\mathbb C}}
\newcommand\RR{\mathbb{R}}
\newcommand\1{\mathbb{1}}
\newcommand\Pb{\mathbf{P}}
\newcommand\Ex{\mathbf{E}}
\newcommand\Liminf{\mathop{\underline{\lim}}\limits}
\newcommand\Limsup{\mathop{\overline{\lim}}\limits}
\def\zs#1{_{\vphantom{T^{T^T}}#1}}
\def\zT{{\vphantom{\widetilde T}T}}

\def\onepicsize{9cm}
\def\onepic#1#2{\nobreak\bigskip\noindent\centerline{%
\hbox to \onepicsize{\vbox{%
\hbox to \onepicsize{\includegraphics*[width=\onepicsize]{#1}}%
\hbox to \onepicsize{\hfill\small #2\hfill}}}}\bigskip}

\begin{document}
\title{Hypotheses Testing: Poisson Versus Self-correcting}

\author{Sergue\"{\i} \textsc{Dachian}\\
{\small Laboratoire de Math\'ematiques}\\
{\small Universit\'e Blaise Pascal, 63177 Aubi\`ere Cedex, France}\\
Yury A. \textsc{Kutoyants}\\
{\small Laboratoire de Statistique et Processus}\\
{\small Universit\'e du Maine, 72085 Le Mans Cedex 9, France}}


\maketitle
\begin{abstract}

We consider the problem of hypotheses testing with the basic simple
hypothesis: observed sequence of points corresponds to stationary Poisson
process with known intensity against a composite one-sided parametric
alternative that this is a self-correcting point process. The underlying
family of measures is locally asymptotically quadratic and we describe the
behavior of score function, likelihood ratio and Wald tests in the asymptotics
of large samples. The results of numerical simulations are presented.
\end{abstract}



\section{Introduction}

The model of self-correcting point process was proposed in 1972 by Isham and
Wescott \cite{IW} to describe a stationary sequence of events $\left\{t_1,
t_2,\ldots \right\}$ which does not have the property of Poisson process of
independence of increments on the disjoint intervals.  To introduce this
processes we denote by $X=\left\{X_t, \;t\geq 0\right\}$ the counting process,
i.e., $X_t$ is equal to the number of events on the time interval
$\left[0,t\right]$. Recall that for a stationary Poisson process with a
constant intensity $S>0$ the increments of $X$ on disjoint intervals are
independent and distributed according to Poisson law
$$
\ds \Pb\left\{X_t-X_s=k\right\}=\frac{S^k\left(t-s\right)^k}{k!}\;
{\rm e}^{-S\left(t-s\right) },\quad 0 \leq s<t,\quad k=0,1,\ldots.
$$
Particularly,
$$
\ds \Pb\left\{X_{t+{\rm d}t}-X_t>0\right\}=S\;{\rm
d}t\;\left(1+o\left(1\right)\right).
$$
For self-correcting point process we have
$$
\ds \Pb\left\{X_{t+{\rm d}t}-X_t>0\; |\; {\cal F}_t\right\}=S\left(t,X_t\right)\;{\rm
d}t\;\left(1+o\left(1\right)\right)
$$
where ${\cal F}_t $ is the $\sigma $-field generated by $\left\{X_s,0\leq
s\leq t\right\}$ and the intensity function
$$
S\left(t,X_t\right)=a\,\psi\left(at -X_t\right),\quad t\geq 0.
$$
Here $a >0$ and the function $\psi\left(\cdot \right)$ satisfies the
following conditions:
\begin{enumerate}
\item $0\leq \psi\left(x\right)<\infty  $ for any $x\in \RR$,
\item there exists a positive constant $c$ such that $ \psi\left(x\right)\geq
c  $ for any $x>0$,
\item $\Liminf_{x\rightarrow \infty } \psi\left(x\right)>1  $,
and $\Limsup_{x\rightarrow -\infty } \psi\left(x\right)<1  $.
\end{enumerate}

Self-correcting processes are called as well stress-release processes (see
\cite{DaVer}, p. 239). This class of processes is widely used as a good
mathematical model for non-poissonian sequences of events. This model was
found especially attractive in the description of earthquakes (see Ogata and
Vere-Jones \cite{VJ-O1}, Lu {\sl at al.} \cite{LHB}).

\bigskip
\noindent
{\bf Example 1.} Let
\begin{equation*}
S\left(t,X_t \right)=\exp \left\{\alpha +\beta \left(t-\varrho
X_t\right)\right\}
\end{equation*}
where $\beta >0, \;\varrho >0$. It is easy to see that the conditions 1--3 are
fulfilled and the point process with such intensity function is
self-correcting.

\bigskip

This model was studied by many authors (see the references in
\cite{DaVer}). Particularly it was shown that under mild conditions there
exists an invariant measure $\mu $ and the law of large numbers (LLN)
\begin{equation}
\label{1}
\frac{1}{T}\int_{0}^{T}h\left(St-X_t\right)\; {\rm d}t\longrightarrow
\int_{}^{}h\left(y\right)\;\mu \left({\rm d}y\right)
\end{equation}
is valid (see Vere-Jones and Ogata \cite{VJ-O2}, Hayashi \cite{Hay}), Zheng \cite{Z}.  Here $h\left(\cdot \right)$
is a continuous, integrable (w.r.t. $\mu $) function and $S>0$ is the rate of the point process.  For the model of
Example 1 we have the LLN if $\rho >0$ and $\beta >0$.

\bigskip

As the self-correcting model is an alternative for the stationary Poisson
process, it is natural and important to test these two hypotheses
 by the
observations $\left\{t_1,t_2, \ldots\right\}$ on the time interval
$\left[0,T\right]$, i.e., to test
$$
S\left(t,X_t\right)=S\;\quad {\rm
versus}\;\quad S\left(t,X_t\right)=a\,\psi\left(at-X_t\right).
 $$
 Remind that the likelihood ratio in this problem has the
following form
\begin{align*}
L\left(X^T\right)=&\exp\left\{\int_{0}^{T}\ln\frac{a\,\psi\left(at-X_{t-}\right)
}{S}\;\left[{\rm d}X_t-S\,{\rm d}t \right]\right.\\
&\left. -\int_{0}^{T}\left[\frac{a\,\psi\left(at-X_{t}\right)}{S}-1-
\ln\frac{a\,\psi\left(at-X_{t}\right) }{S}\;\right]\;S\,{\rm d}t
\right\},
\end{align*}
where $X_{t-} $ is the limit from the left of $X_{t}$ at the point $t$
\cite{LS1}.  Therefore, if the function $a\psi\left(\cdot \right)/S$ is
separated from 1 then the second integral in this representation tends to
infinity and there are many consistent tests. Hence it is more interesting to
compare tests in the situations when the alternatives are {\sl contigous},
i.e. the corresponding sequence of measures are contigous. This corresponds
well to Pitman's approach in hypotheses testing \cite{Pit}. We can have such
situations if $\psi\left(\cdot \right)=S+o\left(1\right) $ with special rates
$o\left(1\right)$. In this work we consider one of such models defined by the
intensity function $S\left(t,X_t\right)=S \psi\left(\vartheta
\left(St-X_t\right)\right)$ where $\vartheta $ is a {\sl small parameter} and
$\psi\left(0\right) =1$. We suppose that the function $\psi\left(\cdot
\right)$ is smooth and we can write
\begin{align*}
&\int_{0}^{T}\left[\psi\left(\vartheta \left(St-X_{t}\right)\right)-1-
\ln\psi\left(\vartheta \left(St-X_{t}\right)\right) \;\right]\;S\,{\rm
d}t=\\
&\qquad \qquad =\frac{\vartheta ^2\dot\psi\left(0\right)^2S}{2}\int_{0}^{T}
\left(St-X_{t}\right)^2{\rm d}t  \;\left(1+o\left(1\right)\right).
\end{align*}
It is easy to see that the rate $\vartheta =\vartheta _T\rightarrow 0$ under
hypothesis $S\left(t,X_t\right)=S  $ is $\vartheta _T\sim T^{-1}$ because
$$
\frac{1}{S\,T^2}\int_{0}^{T}
\left(St-X_{t}\right)^2{\rm d}t =\int_{0}^{1}W_T\left(s\right)^2\,{\rm
d}s\Longrightarrow \int_{0}^{1}W\left(s\right)^2\,{\rm d}s
$$
where $W_T\left(s\right)=\left(ST\right)^{-1/2}\left(S\,Ts-X_{Ts}\right)
\Rightarrow  W\left(s\right)$, and $\left\{W\left(s\right),\;
0\leq s\leq 1\right\}$ is Wiener process.
Note that we put $a=S$, otherwise
\begin{align*}
&\frac{\dot\psi\left(0\right)^2 \vartheta _T^2}{2}\int_{0}^{T}
\left(at-X_t\right)^2\;{\rm d}t=\\
&\quad   =\frac{\dot\psi\left(0\right)^2 \;\vartheta
 _T^2}{2}\int_{0}^{T}\left(\left(a-S\right)t+\sqrt{ST}\;\frac{St-X_t}{\sqrt{ST}}
\right)^2\;{\rm d}t\\
&\quad  =\frac{\dot\psi\left(0\right)^2 \;\vartheta  _T^2}{2}\,T\, \int_{0}^{1}
\left(\left(a-S\right)vT+\sqrt{ST}\;W_T\left(v\right)\right)^2\;{\rm d}v\\
&\quad  =
\frac{\dot\psi\left(0\right)^2 }{6}\;\vartheta  _T^2 \left(a-S
\right)^2\;T^3\;\left(1+o\left(1\right)\right)
 \end{align*}

Therefore, if $a\neq S$, then we have to take $\vartheta _T=uT^{-3/2}$ and to
test the simple hypothesis ${\scr H}_0: u=0$ against ${\scr H}_1: u>0$. In
this case the
family of measures is LAN and the usual construction provides us {\sl
asymptotically uniformly most powerful test} (see, e.g., Roussas
\cite{Rou}). Note that according to \eqref{1} for any fixed
alternative $\vartheta >0$ we have the convergence
$$
\frac{1}{T}\int_{0}^{T}\left(St-X_{t}\right)^2{\rm d}t \longrightarrow
\int_{}^{}y^2\;\mu \left({\rm d}y\right)
$$
which, of course, requires another normalization.

\bigskip

Therefore we consider the problem of hypotheses testing when under hypothesis
${\scr H}_0$ the intensity function is a known constant $S>0$ (Poisson
process) and the alternative ${\scr H}_1$ is one-sided composite:
self-correcting process with intensity function $S\left(t,X_t\right)=S
\psi\left(\vartheta_T \left(St-X_t\right)\right)$, where for convenience of
notation we put $\vartheta _T=u/S\dot\psi\left(0\right)T$ (we suppose that
$\dot\psi\left(0\right)>0 $).  In this case the corresponding likelihood ratio
$Z_T\left(u\right)$ converges to the limit process
$$
Z\left(u\right)=\exp\left\{-u\int_{0}^{1}W\left(s\right)\,{\rm
d}W\left(s\right)-\frac{u^2}{2}\int_{0}^{1}W\left(s\right)^2\,{\rm
d}s \right\},
$$
i.e., the family of measures is {\sl locally asymptotically quadratic}
\cite{LC-Y}. We study three tests: {\sl score function test}, {\sl likelihood
ratio test}, {\sl Wald test} and compare their power functions with the power
function of the {\sl Neyman-Pearson test}.  Note that we calculate all limits
under hypothesis (Poisson process) and we obtain the limit distributions of
the underlying statistics under alternative (self-correcting process) with the
help of Le Cam's Third Lemma. Therefore we do not use directly the conditions
1--3 given above.

The similar limit likelihood ratio process arises in the problem of hypotheses
testing  $u=0$ against  $u>0$ for the time series
$$
X_j=\left(1-\frac{u}{n}\right)\;X_{j-1} + \varepsilon _j,\qquad \quad
j=1,\ldots,n\rightarrow \infty ,
$$
where $\varepsilon _j$ are i.i.d. random variables, $\Ex \varepsilon _j=0, \Ex
\varepsilon _j^2=\sigma ^2$. The asymptotic properties of tests are described
under hypothesis and alternatives by Chan and Wei \cite{C-W} and Phillips
\cite{Phil}. Particularly, the limits of the power functions are given with
the help of Ornstein-Uhlenbeck process
$$
{\rm d}Y_s=-u\,Y_s\;{\rm d}s+{\rm d}W_s,\quad Y_0=0, \qquad 0\leq s\leq 1.
$$
Then Swensen \cite{Swen97} compared these limit powers.

For the model of Example 1 the power functions (for local alternatives) was
studied by Ogata and Vere-Jones \cite{VJ-O1} and by Luschgy \cite{Lus1},
\cite{Lus2}. The limit likelihood ratio and tests are similar to that of the
mentioned above time series problem.  Remind as well that Feigin \cite{Fei}
noted that the same limit likelihood ratio arises in the problem of testing
the simple hypothesis $u=0$ against one-sided alternative $u>0$ by
observations
$$
{\rm d}X_t=-\frac{u}{T}\,X_t\;{\rm d}t+{\rm d}W_t,\quad X_0=0, \qquad 0\leq
t\leq T\rightarrow \infty .
$$

In our case we obtain similar limit expressions for the likelihood ratio and
power functions and compare the errors of tests. The analytical considerations
give us an asymptotic (for large values of $u$) ordering of the tests. The
numerical simulations of the tests show that for the small values of
$\varepsilon $ and for the moderate values of $u$ the power functions of the
likelihood ratio and Wald tests are indistinguishable (from the point of view
of numerical simulations) of the Neyman-Pearson envelope.  This interesting
property was noticed (for $\varepsilon =0.05$) by Eliott {\sl at al.}
\cite{Eli} on the base of $2\cdot 10^3$ simulations. In our work we obtain
similar result having $10^7 $ simulations and we observe for the larger values
of $\varepsilon $ that the asymptotic ordering of the tests holds already
for the moderate values of $u$.

A similar problem of hypotheses testing in the situation, when the
alternative process is self-exciting \cite{Haw} was considered in \cite{DaK}.

\section{Score Function Test}

We observe a trajectory $X^T=\left\{X_t,0\leq t\leq T\right\}$ of a point
process of intensity function $S\left(\cdot ,X_t \right) $ and consider the
problem of testing the simple hypothesis against close one sided composite
alternative
\begin{eqnarray}
\label{3}
&&{\scr H}_0:\quad\qquad S\left(t,X_t \right)=S_*,\\
\label{4}
&&{\scr H}_1:\quad\qquad S\left(t,X_t \right)=S_*\; \psi\left(\vartheta
_T\left[S_*t-X_t\right]\right), \quad \vartheta _T>0,
\end{eqnarray}
where $\vartheta _T$ is a small parameter, the value $S_*$ and the function
$\psi \left(\cdot \right) $ are known.
The problem is regular in the following sense.

\bigskip

{\bf Condition} ${ \cal A}.$ {\sl The function $\psi\left(x \right), x\in \RR$
is positive,
continuously differentiable  at the point $x=0$, $\psi\left(0\right)=1$ and
$\dot \psi\left(0\right)>0$.}

\bigskip

 The rate of convergence $\vartheta _T\rightarrow 0$ is
chosen such that the likelihood ratio $L\left(\vartheta _T, X^T\right) $ is asymptotically non degenerate.  In the
case $\dot \psi\left(0\right)<0$ we need to change just one sign in the test.  This leads us to the
reparametrization
$$
\vartheta _T=\frac{u}{S_*\;\dot \psi\left(0\right)\;T}, \;u\geq 0
$$
 and to the corresponding  hypotheses testing problem
\begin{eqnarray}
\label{5}
&&{\scr H}_0:\quad\qquad u=0,\\
\label{6}
&&{\scr H}_1:\quad\qquad u>0.
\end{eqnarray}
 Therefore, we observe a Poisson process of intensity $S_*$ under hypothesis
 ${\scr H}_0  $ and the point process under alternative ${\scr H}_1 $ has intensity function
$$
 S\left(t,X_t
 \right)=S_*+\frac{u}{T}\left(S_*t-X_{t}\right)+o\left(T^{-1/2}\right).
$$

Let us fix $\varepsilon \in \left(0,1\right)$ and denote by ${\scr
K}_\varepsilon $ the class of test functions $\phi_T\left(X^T\right)$ of
asymptotic size $\varepsilon $, i.e., for $\phi_T\in {\scr K}_\varepsilon $
we have
$$
\lim_{T\rightarrow \infty }\Ex_0\;\phi_T\left(X^T\right)=\varepsilon.
$$
As usual, $ \phi_T\left(X^T\right)$ is the probability to accept the
hypothesis ${\scr H}_1 $ having observations $X^T$. The corresponding power
function is
$$
\beta _\zT\left(u,\phi_\zT\right)=\Ex_u\;\phi_T\left(X^T\right),\quad u\geq 0.
$$

Let us introduce the statistic
\begin{align}
\Delta
_T\left(X^T\right)&=\frac{1}{S_*\;T}
\int_{0}^{T}\left(S_*t-X_{t-}\right)\;\left[{\rm d}X_t-S_*\;{\rm d}t\right]\nonumber\\
&=\frac{X_T-(X_T-S_*T)^2}{2\,S_*T}.
\label{7}
\end{align}
This equality follows from the elementary representation  (see,
e.g. \cite{Kut84}, Lemma 4.2.1) for the centered Poisson process $\pi _t=X_t-S_*\;t$
\begin{align*}
\pi _T^2=2\int_{0}^{T}\pi _{t-}\;{\rm d}\pi _t+\pi _T+S_*T.
\end{align*}
which obviously is equivalent to
$$
\frac{1}{T}\int_{0}^{T}\pi _{t-}\;{\rm d}\pi _t=\frac{\pi
_T^2-X_T}{2T}.
$$

Define as well  two random variables
$$
\Delta(W)
=\frac{1}{2}\left(1-W\left(1\right)^2\right)=-\int_{0}^{1}W\left(s\right){\rm
d} W\left(s\right),\qquad {\rm J}(W)
=\int_{0}^{1}W\left(s\right)^2{\rm d}s ,
$$
where $\left\{W\left(s\right), 0\leq s\leq 1 \right\}$  is standard Wiener
process.

Remind that the likelihood ratio in this problem has the following form
\cite{LS1} $\left(\gamma =S_*\dot\psi\left(0\right)\right)$
\begin{align}
\nonumber
&L\left(\frac{u}{\gamma T},X^T\right)=\exp\left\{\int_{0}^{T}\ln
\psi\left(\frac{u}{\gamma T}\left(S_*\;t-X_{t-} \right)\right)\;
\left[{\rm d}X_t-S_*\;{\rm d}t\right]\right.\\
&\quad \left.-\int_{0}^{T} \left[\psi\left(\frac{u}{\gamma T}\left(S_*\;t-X_t
\right)\right)-1-\ln \psi\left(\frac{u}{\gamma T}\left(S_*\;t-X_t
\right)\right)    \right]S_*\;{\rm d}t\right\}.
\label{8}
\end{align}

Therefore the direct differentiation w.r.t. $u$ at the point $u=0$ gives us the
introduced above statistic
$$
\left.\frac{\partial }{\partial\; u}\ln L\left(\frac{u}{\gamma T},X^T\right)\right|_{u=0}=\Delta
_T\left(X^T\right) .
$$
Below we denote
$$
a_\varepsilon=\frac{1-z^2_{\frac{1-\varepsilon}{2}}}{2}\qquad {\rm and}\qquad
h\left(u\right)=\sqrt{\frac{2u}{1-e^{-2u}}},
$$
where $z_a$ is $1-a$ quantile of standard Gaussian law, i.e., $\Pb\left(\zeta
>z_a\right) =a$, for $\zeta \sim {\cal N}\left(0,1\right)$.

We have the following result.

\begin{theorem}
\label{T1}
Let the Condition ${\cal A}$ be fulfilled, then the score function test
\begin{equation}
\label{9}
\phi_T^*\left(X^T\right)=\1\zs{\left\{\Delta
_T\left(X^T\right)>a_\varepsilon \right\}}
\end{equation}
 belongs to the class ${\scr
K}_\varepsilon $ and  for any $u_*>0$ its power function
\begin{align}
\label{10}
\beta _T(u_*,\phi_T^*)\rightarrow \beta^*(u_*)=\Pb\left\{\left|\zeta
\right|\leq  h\left(u_*\right)\;z_{\frac{1-\varepsilon
}{2}} \right\} .
\end{align}

\end{theorem}

\bigskip
\noindent
{\bf Proof.}
Under hypothesis ${\scr H}_0$ the value  $X_T$ is a poissonian random variable
with parameter $S_*T$. Therefore we have
immediately
$$
\frac{X_T}{S_*T}\longrightarrow 1,\qquad \qquad \frac{X_T-ST}{\sqrt{S_*T}
}\Longrightarrow W\left(1\right)\sim {\cal N}\left(0,1\right)
$$
and
$
\Delta _T\left(X^T\right)\Longrightarrow \Delta \left(W\right)
$
as $T\rightarrow \infty $. Hence
$$
\Pb_0\left\{\Delta _T\left(X^T\right)>a_\varepsilon \right\}\longrightarrow
\Pb\left\{\Delta \left(W\right)> \frac{1-z^2_{\frac{1-\varepsilon}{2}}}{2} \right\} =
\Pb\left\{\left|\zeta \right| < z_{\frac{1-\varepsilon}{2}}\right\}=\varepsilon .
$$
This provides  $ \phi_T^*\in {\scr K}_\varepsilon $.

\bigskip

To study the power $\beta _T(u_*,\phi_T^*)$ we would like to use the Third
Le Cam Lemma \cite{LC-Y}, \cite{Str}.
Therefore we need first to show the joint weak convergence
\begin{align}
\label{11}
{\cal L}_0\left(\Delta _T , l _T\left(u\right) \right)\Longrightarrow  {\cal
L}\left(\Delta(W),u\Delta(W)-\frac{u^2}{2}{\rm J}(W)\right)
\end{align}
where $ l _T\left(u\right)=\ln L\left(\frac{u}{\gamma
T},X^T\right)$.

To verify \eqref{11} we denote
$$
l _T^*\left(u\right)=u\;\Delta_T\left(X^T\right)-\frac{u^2}{2}\; J_T\left(X^T\right),
$$
where
$$
J_T\left(X^T\right)=\frac{1}{S_*\;T^2}\int_{0}^{T}\left(S_*
t-X_t\right)^2\;{\rm d}t
$$
and show that
\begin{align}
\label{12}
{\cal L}_0\left( l _T^*\left(u\right) \right)\Longrightarrow  {\cal
L}\left(u\Delta(W)-\frac{u^2}{2}{\rm J}(W)\right).
\end{align}
Then \eqref{11} will follow from the convergence
\begin{equation}
\label{13}
l _T^*\left(u_T\right)-l _T\left(u_T\right)\rightarrow 0
\end{equation}
for any bounded sequence $u_T$.

\bigskip

\begin{lemma}
\label{L1}
 \begin{align}
\label{14}
{\cal L}_0\left\{\Delta
_T\left(X^T\right),\; J_T\left(X^T\right)\right\} \Longrightarrow
\left(-\int_{0}^{1}W\left(s\right)\;{\rm d}W\left(s\right),
\int_{0}^{1}W\left(s\right)^2\;{\rm d}s\right).
\end{align}
\end{lemma}

\bigskip

\noindent
{\bf Proof.} Let us put  $W
_T\left(s\right)=\left(S_*T\right)^{-1/2}\pi _{sT},\;s\in
\left[0,1\right]$. Then
$$
\Ex_0W_T\left(s\right)=0,\quad \Ex_0\left[W _T\left(s_1\right)W
_T\left(s_2\right)\right]= \min\left(s_1,s_2\right)
$$
and we have
$$
J_T\left(X^T\right)=\frac{1}{S_*\;T^2}\int_{0}^{T}\pi
_t^2\;{\rm d}t =\int_{0}^{1}W
_T\left(s\right)^2\;{\rm d}s.
$$
Using the standard arguments we verify (well-known fact)
 that for any collection $\left\{s_1,\ldots, s_k\right\}$ we have the weak
 convergence (as $T\rightarrow \infty $) of the vectors
$$
\Bigl(W _T\left(s_1\right),\ldots,W
_T\left(s_k\right)\Bigr)\Longrightarrow
\Bigl(W\left(s_1\right),\ldots,W\left(s_k\right)\Bigr) .
$$
Moreover the following estimate holds
\begin{align*}
&\left(\Ex_0\left|W _T\left(s_1\right)^2-W
 _T\left(s_2\right)^2\right|\right)^2\leq \\
&\qquad \qquad \leq
 \Ex_0\left|W_T\left(s_1\right)-W _T\left(s_2\right)\right|^2
 \Ex_0\left|W _T\left(s_1\right)+W _T\left(s_2\right)\right|^2\leq
 4\,\left|s_2-s_1\right| .
\end{align*}

Hence (see Gikhman and Skorohod \cite{GS}, Section IX.7) we have the
convergence (in distribution) of integrals
$$
\int_{0}^{1}W_T\left(s\right)^2\;{\rm d}s\Longrightarrow
\int_{0}^{1}W\left(s\right)^2\;{\rm d}s
$$
and
$$
\Delta
_T\left(X^T\right)=\frac{1-W_T\left(1\right)^2}{2}\;
\left(1+o\left(1\right)\right) \Longrightarrow
\frac{1-W\left(1\right)^2}{2}=-\int_{0}^{1}W\left(s\right){\rm
d}W\left(s\right) .
$$
It is easy to see that we have the same time the joint convergence too because
from the given above proof it follows that for
any $\lambda _1,\lambda _2$
$$
\lambda _1 W_T\left(1\right)^2+\lambda _2\int_{0}^{1}W_T\left(s\right)^2\;{\rm
d}s\Longrightarrow
\lambda _1 W\left(1\right)^2+\lambda _2\int_{0}^{1}W\left(s\right)^2\;{\rm d}s .
$$
Therefore the Lemma \ref{L1} is proved.

\bigskip

Our goal now is to establish a slightly more strong than \eqref{13}   relation
\begin{align}
\label{15}
l_T\left(u_T\right)=u_T\;\Delta _T\left(X^T\right)\;\left(1+o\left(1\right)\right)
-\frac{u_T^2}{2}\int_{0}^{1}W_T\left(s\right)^2\;{\rm
d} s\;\left(1+o\left(1\right)\right)
\end{align}
where $o\left(1\right)\rightarrow 0$ for any sequence $u_T\in \UU_T$ with
$\UU_T=\{u:\; 0\leq
u<\frac{\sqrt{S_*T}}{\ln T}\}$.

We can write
\begin{align*}
&l_T^*\left(u\right)-l_T\left(u\right)=\int_{0}^{T}
\left[-\frac{u\;W_T\left(\frac{t}{T}\right)}{\sqrt{S_*T}}- \ln
\psi\left(\frac{-uW_T\left(\frac{t}{T}\right)}{\dot\psi\left(0\right)\sqrt{S_*T}}\right)
\right]\;{\rm
d}\pi _t\\
&\quad  -\int_{0}^{T}
\left[\frac{u^2\;W_T\left(\frac{t}{T}\right)^2}{2S_*\;T}-\psi\left(
\frac{-uW_T\left(\frac{t}{T}\right)}{\dot\psi\left(0\right)\sqrt{S_*T}}
\right)+1+\ln \psi\left(\frac{-uW_T\left(\frac{t}{T}\right)}{
\dot\psi\left(0\right)\sqrt{S_*T}}\right)   \right]S_*\;{\rm d}t\\
&\quad\equiv u\; \delta _{1,T}-\frac{u^2}{2}\;\delta _{2,T}
\end{align*}
with obvious notation. Remind that $u>0$. Using Lenglart inequality we obtain for the first term
\begin{align*}
&\Pb_0\left\{ \left|\delta _{1,T}\right|>a\right\}\leq
\frac{b}{a}\\
& \qquad +\Pb_0\left\{\int_{0}^{1} \left[W_T\left(s\right)+\frac{\sqrt{S_*T}}{u}\ln
\psi\left(\frac{-uW_T\left(s\right)}{\dot\psi\left(0\right)\sqrt{S_*T}}\right)
\right]^2\;{\rm d}s>b\right\}
\end{align*}
for any $a>0$ and $b>0$. Now expanding the functions $\psi\left(\cdot \right)$
 we obtain
$$
\psi\left(\frac{-uW_T\left(s\right)}{\dot\psi\left(0\right)\sqrt{S_*T}}\right)=
1-\frac{uW_T\left(s\right)}{\dot\psi\left(0\right)\sqrt{S_*T}}\;\dot\psi\left(\frac{-\tilde
uW_T\left(s\right)}{\dot\psi\left(0\right)\sqrt{S_*T}}\right)
$$
where $\tilde u\leq u$. Introduce the set
$$
\CC_T=\left\{\omega :\qquad \sup_{0\leq s\leq 1}\left|W_T\left(s\right)
\right|\leq \dot\psi\left(0\right)\sqrt{\ln T}\right\}
$$
and note that for $\omega \in \CC_T$ we have the estimate
$$
\sup_{u\in \UU_T}\sup_{0\leq s\leq 1} \frac{ u\left|W_T\left(s\right)
\right|}{ \dot\psi\left(0\right)\sqrt{S_*T} }\leq \;\frac{1}{\sqrt{\ln T}}
$$
Hence for all $u\in \UU_T$  on this set  we can write
$$
\sup_{0\leq s\leq 1 }\left|\dot\psi\left(0\right)-\dot\psi\left(\frac{-\tilde
uW_T\left(s\right)}{\dot\psi\left(0\right)\sqrt{S_*T}}\right)
\right|\leq\sup_{\left|v\right|\leq \left(\ln T\right)^{-1/2}
}\left|\dot\psi\left(0\right)-\dot\psi\left(v\right) \right| = h_T\rightarrow 0
$$
as $T\rightarrow \infty $ because the derivative is continuous at the point $v=0$.

Let us denote
$u_s=\frac{uW_T\left(s\right)}{\dot\psi\left(0\right)\sqrt{S_*T}}$.
 Using the expansion of the  logarithm
$$
\ln \left(\psi\left(-u_s\right)\right)=\ln
\left(1-u_s\dot\psi\left(-\tilde
u_s\right)\right)=-\frac{u_s\dot\psi\left(-\tilde u_s\right)}{1-\tilde{\tilde
{u}}_s\dot\psi\left(-\tilde u_s\right)} .
$$
  we obtain the following  estimate
\begin{align*}
&\Pb_0\left\{\int_{0}^{1}
\left[W_T\left(s\right)+\frac{\sqrt{S_*T}}{u}\ln
\psi\left(-u_s  \right)
\right]^2\;{\rm d}s>b\right\} \leq  \Pb_0\left\{\CC_T^c\right\}+\\
&\quad  +
\Pb_0\left\{\int_{0}^{1}W_T\left(s\right)^2\left(1-\frac{\dot\psi\left(-\tilde u_s\right) }{\dot\psi\left(0\right)\;\left(1-\tilde{\tilde
{u}}_s\dot\psi\left(-\tilde u_s\right)\right)}\right)^2{\rm
d}s>b, \CC_T\right\}.
\end{align*}
Remind that $W_T\left(s\right)$ is martingale, hence by Doob inequality we
have
$$
\Pb_0\left\{\CC_T^c\right\}\leq
\Pb_0\left\{ \left|W_T\left(1\right)\right|> \dot\psi\left(0\right)\sqrt{\ln T}
\right\}\leq \frac{1 }{\dot\psi\left(0\right)^2\ln T}.
$$
For the second probability after elementary estimates we obtain
\begin{align*}
&\Pb_0\left\{\int_{0}^{1}W_T\left(s\right)^2\left(1-\frac{\dot\psi\left(-\tilde
u_s\right) }{\dot\psi\left(0\right)\;\left(1-\tilde{\tilde
{u}}_s\dot\psi\left(-\tilde u_s\right)\right)}\right)^2{\rm
d}s>b, \CC_T\right\}\leq \\
&\quad\quad  \leq \Pb_0\left\{C\,\int_{0}^{1}W_T\left(s\right)^2{\rm
d}s\;\left(h_T^2+\frac{1}{\ln T}\right)>b\right\}\leq \frac{C}{2b}\,\left(h_T^2+\frac{1}{\ln
T}\right)
\end{align*}
with some constant $C>0$.
Recall that by Tchebyshev inequality
$$
\Pb_0\left\{\int_{0}^{1}W_T\left(s\right)^2{\rm
d}s>A\right\}\leq \; \frac{1}{2A}.
$$

Therefore, if we take $b=a^2$  then for any $a>0$
$$
\Pb_0\left\{ \left|\delta _{1,T}\right|>a\right\}\longrightarrow 0
$$
as $T\rightarrow \infty $.

The similar arguments allow to prove the convergence
$$
\Pb_0\left\{ \left|\delta _{2,T}\right|>a\right\}\longrightarrow 0
$$
too.

\bigskip

Therefore, the likelihood ratio $Z_T\left(u\right)=L\left(\frac{u}{\gamma T},X^T\right),
u\geq 0$ is (under hypothesis ${\scr H}_0$) {\sl
locally asymptotically quadratic} (LAQ) \cite{LC-Y}, because
\begin{align}
\label{16}
Z_T\left(u\right)\Longrightarrow Z\left(u\right)=\exp\left\{-u
\int_{0}^{1}W\left(s\right)\;{\rm d}W\left(s\right)
-\frac{u^2}{2}\int_{0}^{1}W\left(s\right)^2\;{\rm d}s\right\}.
\end{align}
Moreover, we have the convergence
$l_T^*\left(u_T\right)-l_T\left(u_T\right)\rightarrow 0$ for any bounded
sequence of $u_T\in \UU_T$.
 Note that the random function $Z\left(u\right) $
is the likelihood ratio in the hypotheses testing  problem
\begin{eqnarray*}
&&{\scr H}_0:\quad\qquad u=0,\\
&&{\scr H}_1:\quad\qquad u>0,
\end{eqnarray*}
by observations of Ornstein-Uhlenbeck process
\begin{equation}
\label{17}
{\rm d}Y\left(s\right)=-uY\left(s\right)\;{\rm d}s+{\rm
d}W\left(s\right),\quad Y\left(0\right)=0,\qquad 0\leq s\leq 1
\end{equation}
under hypothesis $u=0$.

This limit for the likelihood ratio under alternative can be obtained directly
as follows. Let us denote
$$
Y_T\left(s\right)=
\frac{X_{sT}-sS_*T}{\sqrt{S_*T} },\qquad 0\leq s\leq 1.
$$
Then using the representation
$$
X_t=S\,\int_{0}^{t}\psi\left(\vartheta _T\left[S_* r-X_r\right]\right)\;{\rm d}r+M_t
$$
where $M_t$ is local martingale
  and expansion of  the function $\psi\left(\cdot \right)$ at the
vicinity of $0$ we obtain the equation
$$
Y_T\left(s\right)=-u\int_{0}^{s}\frac{\dot\psi\left(g_v\right)}{\dot
\psi\left(0\right)}Y_T\left(v\right)\;{\rm
d}v + V_T\left(s\right),\quad Y_T\left(0\right)=0,\qquad 0\leq s\leq 1
$$
where $V_T\left(s\right)$ is local martingale and $g_v= \frac{-\tilde
u}{\dot\psi\left(0\right)\sqrt{S_*T}} Y_T\left(v \right)\rightarrow 0$. The
central limit theorem for local martingales provides the convergence
$V_T\left(s\right)\Longrightarrow W\left(s\right)$. Hence the process
\eqref{17} is the limit (in distribution) of $Y_T\left(s\right)$. Moreover
from \eqref{8} we have
$$
\Delta
_T\left(X^T\right)=\frac{Y_T\left(1\right)}{2\sqrt{S_*T}}+
\frac{1-Y_T\left(1\right)^2}{2}\Longrightarrow \frac{1-Y\left(1\right)^2}{2}.
$$
This limit of the statistic $\Delta
_T\left(X^T\right)$ follows from the Third Le Cam Lemma as well. Particularly,
for any continuous bounded function $H\left(\cdot \right)$
\begin{align*}
&\Ex_uH\left(\Delta _T\left(X^T\right)\right)=\Ex_0
\left[Z_T\left(u\right)H\left(\Delta _T\left(X^T\right)\right)\right]
\longrightarrow \\
&\qquad \longrightarrow \Ex_0
\left[Z\left(u\right)H\left(\Delta \left(W\right)\right)\right]
=\Ex_uH\left(\Delta \left(Y\right)\right),
\end{align*}
where
$$
\Delta \left(Y\right)=-\int_{0}^{1}Y\left(s\right)\;{\rm
d}Y\left(s\right)=\frac{1-Y\left(1\right)^2}{2}.
$$

Hence under alternative $\left(\vartheta
_T=u_*/\gamma T\right) $ we have the convergence
$$
\beta_T \left(u_*,\phi_T^*\right)\longrightarrow
\Pb_{u_*}\left\{\left|Y\left(1\right)\right|\leq z_{\frac{1-\varepsilon }{2}}
\right\}= \Pb\left\{\left|W\left(1\right)\right|\leq z_{\frac{1-\varepsilon
}{2}}\sqrt{\frac{2u_*}{1-e^{-2u_*}}} \right\}
$$
because
$$
Y\left(1\right)=\int_{0}^{1}e^{-u\left(1-s\right)}\;{\rm d}W\left(s\right)\sim
{\cal N}\left(0,\frac{1-e^{-2u_*}}{2u_*} \right)
$$

This proves \eqref{10}.

\bigskip

Theorem~\ref{T1} is {\sl asymptotic in nature}, and it is interesting to see
the powers of the score function test for the moderate values of $T$ and
especially to compare them with the limit power functions. This can be done
using numerical simulations.

We consider the model of Example 1 with $S_*=1$ and $\psi(t)=e^t$. This yields
the intensity function
$$
S\left(u,t,X_t \right)=\exp\left(\frac uT\left[t-X_t\right]\right), \qquad
u\geq 0, \quad 0\leq t\leq T.
$$

In Figure~1 we represent the power function of the score function test
$\phi^*_\zT$ of asymptotic size $0.05$ given by
$$
\beta_\zT\left(u,\phi^*_\zT\right)=\Pb_u\left\{\Delta_T
\left(X^T\right)>a_{0.05}\right\},\qquad 0\leq u\leq 20,
$$
for $T=100$, $300$ and $1000$, as well as the limiting power function
$\beta^*(\cdot)$ given by the formula~\eqref{10}.

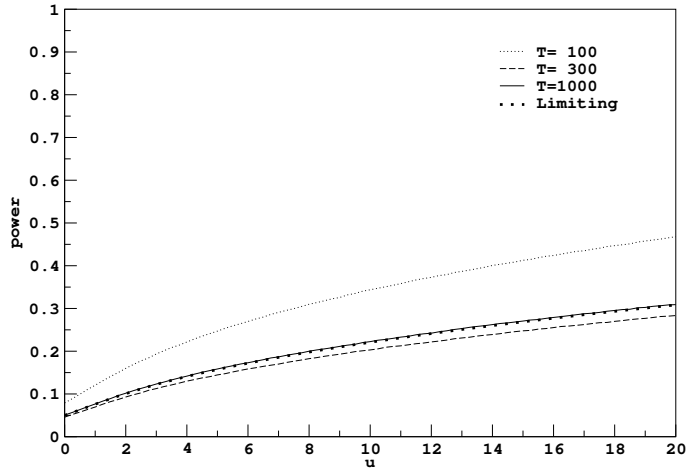
\onepic{SCpowNB.eps}{Fig.~1: Power of the score function test}

The function $\beta_\zT(\cdot,\phi^*_\zT)$ is estimated in the following way.
We simulate (for each value of $u$) $M=10^6$ trajectories $X^T_j$,
$j=1,\ldots,M$ of self-correcting process of intensity
$S\left(u,t,X_t\right)$ and calculate $\Delta_j=\Delta_T(X^T_j)$. Then we
calculate the empirical frequency of accepting the alternative hypothesis
$$
\frac1M\sum_{j=1}^{M}\1\zs{\left\{\Delta_j>a_{0.05}
\right\}}\approx\beta_\zT(u,\phi^*_\zT).
$$
Note that for $T=1000$ the limiting power function is practically
attained. Note also that for $T=100$ the size of the test is $0.079$ which
explains the position of the corresponding curve.

Remind that score-function test is locally optimal \cite{Cap}.

\section{The Likelihood Ratio Test  and the Wald Test}

Let us study two other well-known tests: the {\sl likelihood ratio test}
$\bar\phi_T$ based on the maximum of the likelihood ratio function and the
{\sl Wald test} $\hat\phi_T$ based on the MLE $\hat\vartheta _\zT$.

Remind that the log-likelihood ratio formula is
\begin{align*}
\ln L\left(\vartheta ,X^T\right)&=\int_{0}^{T}\ln \psi \left(\vartheta
\left(S_*t-X_{t-}\right)\right) \;\left[{\rm d}X_t-S_*\,{\rm d}t\right]\\
&\quad -\int_{0}^{T}\left[\psi \left(\vartheta
\left(S_*t-X_{t-}\right)\right)-1- \ln \psi \left(\vartheta
\left(S_*t-X_{t-}\right)\right)\right]\;S_*\,{\rm d}t
\end{align*}
and the likelihood ratio test is based on the statistic
$$
\delta _T\left(X^T\right)=\sup_{\vartheta\in \Theta }L\left(\vartheta ,X^T\right),
$$
where $\Theta $ is the set of values of $\vartheta $ under alternative.
 The test is given by the decision function
$$
\bar\phi_T\left(X^T\right)=\1\zs{\left\{\delta _T\left(X^T\right)>\tilde
b_\varepsilon \right\}}
$$
where the threshold $\tilde b_\varepsilon $ is chosen from the condition
$\bar\phi_T\in {\scr K}_\varepsilon $.

Note that $\delta _T\left(X^T\right)= L\left(\hat\vartheta_T ,X^T\right)$ as well,
 where $\hat\vartheta_T $ is the maximum likelihood estimator of the parameter
 $\vartheta $.

The reparametrization $\vartheta =\vartheta _T=u/\gamma T$ reduces the problem
\eqref{3}-\eqref{4} to \eqref{5}-\eqref{6} and we have to precise the region
of {\sl local alternatives}.  In the traditional approach of {\sl locally
asymptotically uniformly most powerful tests} \cite{Rou} (regular case) to
check the optimality of a test $\phi_\zT$ we compare the power function $\beta
_T\left(u,\phi_\zT\right)$ with the power function of the Neyman-Pearson test
on the compacts $0\leq u\leq K$ for any $K>0$. For these values of $u$ the
alternatives are always {\sl contigous}. To consider the similar class of
alternatives in our case is not reasonable because the constant $\tilde
b_\varepsilon $ became dependent of $K$. Indeed if we take the test function
$$
\bar\phi_\zT\left(X^T\right)=\1\zs{\Bigl\{\sup\limits_{0<u\leq K}
Z_T\left(u\right)>\tilde b_\varepsilon\Bigr\}} ,\qquad \qquad
Z_T\left(u\right)=L\left(\frac{u}{\gamma T} ,X^T\right),
$$
 then the
condition   $\bar\phi_T\in {\scr K}_\varepsilon $ implies
$\tilde b_\varepsilon=\tilde b_\varepsilon\left(K\right)$.
 Therefore we suppose that
$K=K_T=\frac{\sqrt{S_*T}}{\ln T}\rightarrow \infty $.

Finally, we have the following hypotheses testing problem
\begin{eqnarray}
\label{18}
&&{\scr H}_0:\quad\qquad u=0,\\
\label{19}
&&{\scr H}_1:\quad\qquad u=u_*\in \UU_T
\end{eqnarray}

Therefore, to study
$$
\bar\phi_T\left(X^T\right)=\1\zs{\biggl\{\sup\limits_{u\in \UU_T
}Z_T\left(u\right)>\tilde b_\varepsilon
\biggr\}}
$$
we need to describe the asymptotics of its
errors  under hypothesis ${\scr H}_0$ and alternatives  ${\scr
H}_1$ with $\vartheta =\frac{u_*}{\gamma T},\; u_*\in \UU_T$.

Below
$$
 \Lambda(W)= \frac{\Delta(W)}{\sqrt{2{\rm J}(W) } }.
$$

\begin{theorem}
\label{T2} Let us suppose that  condition ${\cal A}$ is fulfilled and the
value $b_\varepsilon $ is solution of the equation
\begin{equation}
\label{thrB}
\Pb\left(\Lambda(W) > b_\varepsilon\right)=\varepsilon .
\end{equation}
Then the test $\bar\phi_T $  with $\tilde b_\varepsilon=e^{b_\varepsilon^2} $ belongs to $
{\scr K}_\varepsilon $ and its power function converges to the following limit
\begin{align*}
\beta \left(u_*,\bar\phi_T\right)\longrightarrow \hat\beta\left(u_*\right)
=\Pb\left\{ \Lambda(Y_{u_*} ) > b_\varepsilon \right\},
\end{align*}
where
$$
 \Lambda(Y_{u_*})= \frac{\Delta(Y_{u_*})}{\sqrt{2{\rm J}(Y_{u_*}) } }=
\frac{1-Y_{u_*}\left(1\right)^2}{\sqrt{8\;{\rm J}(Y_{u_*}) } }.
$$
and $Y_{u_*}=\left\{Y_{u_*}\left(s\right),0\leq s\leq 1\right\}$ is Ornstein-Uhlenbeck
process \eqref{17} with $u=u_*$.
\end{theorem}
\noindent {\bf Proof.} The log-likelihood process $l_T\left(u\right)=\ln
Z_T\left(u\right)$ admits (under hypothesis ${\scr H}_0 $) the representation \eqref{15}
\begin{align}
\label{20}
l_T\left(u\right)=u\;\Delta _T\left(X^T\right)\;\left(1+\delta
_{1,T} \right)-\frac{u^2}{2}\;{\rm J}_T\left(X^T\right)\;\left(1+\delta_{2,T} \right)
\end{align}
where $\delta_{i,T}\rightarrow 0$ uniformly on $u\in \UU_T$. Hence
\begin{align*}
\Lambda _T\left(X^T\right)^2\equiv\sup_{u\in
\UU_T}l_T\left(u\right)\Longrightarrow \frac{\Delta
\left(W\right)^2}{2{\rm J}\left(W\right)}
\end{align*}
and we have
\begin{align*}
\Ex_0\bar\phi_T\left(X^T\right)=\Pb_0\left\{ \sup_{u\in
\UU_T}l_T\left(u\right)>b_\varepsilon^2 \right\}\longrightarrow
\Pb\left(\Lambda(W)  >b_\varepsilon \right)=\varepsilon .
\end{align*}

Let us fix an alternative $u=u_*$. We have the convergence
\begin{align}
\label{21}
{\cal L}_0\left\{\Lambda _T\left(X^T\right),l_T\left(u_*\right)
\right\}\Longrightarrow {\cal
L}\left\{\Lambda(W) , u_*\;\Delta(W)-\frac{u_*^2}{2}\;{\rm
J}(W)\right\} .
\end{align}

The convergence \eqref{21} allows us to apply Third Le Cam's Lemma as
follows: for any bounded continuous
function  $H\left(\cdot \right)$
\begin{align*}
&\Ex_{u_*}H\left(\Lambda  _T\left(X^T\right)\right)=\Ex_0
\left[Z_T\left({u_*}\right)H\left(\Lambda  _T\left(X^T\right)\right)\right]
\longrightarrow \\
&\qquad \longrightarrow \Ex_0
\left[Z\left({u_*}\right)H\left(\Lambda \left(W\right)\right)\right]
=\Ex_{u_*}H\left(\Lambda \left(Y_{u_*}\right)\right).
\end{align*}

Hence
\begin{align*}
&\beta\left(u_*,\bar\phi_T\right)=\Pb_{u_*}\left\{\sup_{u\in
\UU_T}l_T\left(u\right)>b_\varepsilon^2 \right\}\longrightarrow
\Pb_{u_*}\left\{\Lambda \left(Y_{u_*}\right) >b_\varepsilon  \right\}.
\end{align*}

This completes the proof of the theorem \ref{T2}.

\bigskip

Let us note, that the threshold $b_\varepsilon$ is given implicitly as the
solution of the equation~\eqref{thrB}.  In the following table we give some
values of $b_\varepsilon$ obtained using numerical simulations.

\begin{center}
\begin{tabular}{|l|c|c|c|c|c|c|}
\hline
$\varepsilon$&0.01&0.02&0.03&0.04&0.05&0.1
\\\hline
$b_\varepsilon$&1.814&1.636&1.524&1.440&1.373&1.144
\\\hline
\end{tabular}
\end{center}

These thresholds are obtained by simulating $M=10^7$ trajectories on $[0,1]$
of a standard Wiener process, calculating for each of them the quantity
$\Lambda(W)$ and taking $(1-\varepsilon)M$-th greatest between them.

\bigskip

The next test usually studied in such hypotheses testing problems is the Wald
test
$$
\hat\phi_\zT\left(X^T\right)=\1\zs{\left\{\gamma T \hat\vartheta _T\geq
c_\varepsilon \right\}}
$$
where $\hat\vartheta_T$ is the maximum likelihood estimator of $\vartheta$.

Below
$$
\Gamma(W)= \frac{\Delta(W)}{{\rm J}(W)}.
$$

\begin{theorem}
\label{T3} Let us suppose that  condition ${\cal A}$ is fulfilled and the
value $c_\varepsilon $ is solution of the equation
\begin{equation}
\label{thrC}
\Pb\left(\Gamma (W) >c_\varepsilon \right)=\varepsilon .
\end{equation}
Then the test $\hat\phi_T $ belongs to $ {\scr K}_\varepsilon $ and its power
function for any alternative $u_*$ converges to the following limit
\begin{align*}
\beta \left(u_*,\hat\phi_T\right)\longrightarrow \hat\beta\left(u_*\right)
=\Pb\left\{ \Gamma(Y_{u_*}) >c_\varepsilon \right\},
\end{align*}
where
$$
\Gamma(Y_{u_*})= \frac{\Delta(Y_{u_*})}{{\rm J}(Y_{u_*})}=-u_*+
\frac{\int_{0}^{1} Y_{u_*}(s) \,{\rm d}W(s)}{{\rm J}(Y_{u_*})}\,.
$$
and $Y_{u_*}$ is the same as in Theorem
\ref{T2}.
\end{theorem}
\noindent {\bf Proof.} The proof follows immediately from the representation
\eqref{20}, because
\begin{align*}
&\Pb_0^{\left(T\right)}\left\{\gamma T \hat\vartheta _T\geq c_\varepsilon
\right\}=\Pb_0^{\left(T\right)}\left\{\sup_{0\leq u\leq c_\varepsilon }
Z_\zT\left(u\right)<\sup_{ u> c_\varepsilon,u\in \UU_T }Z_\zT\left(u\right)
\right\}\longrightarrow \\
&\qquad  \longrightarrow \Pb_0\left\{\sup_{0\leq u\leq c_\varepsilon }
Z\left(u\right)<\sup_{ u> c_\varepsilon }Z\left(u\right)\right\}=
\Pb\left\{\Gamma (W) >c_\varepsilon \right\}=\varepsilon
\end{align*}
and (under alternative $u=u_*$)
\begin{align*}
&\Pb_{u_*}^{\left(T\right)}\left\{\gamma T \hat\vartheta _T\geq c_\varepsilon
\right\}=\Pb_{u_*}^{\left(T\right)}\left\{\sup_{0\leq u\leq c_\varepsilon }
Z_\zT\left(u\right)<\sup_{ u> c_\varepsilon,u\in \UU_T }Z_\zT\left(u\right)
\right\}\longrightarrow \\
&\qquad  \longrightarrow \Pb_{u_*}\left\{\sup_{0\leq u\leq c_\varepsilon }
Z\left(u\right)<\sup_{ u> c_\varepsilon }Z\left(u\right)\right\}=
\Pb\left\{\Gamma(Y_{u_*})>c_\varepsilon \right\}=
\hat\beta\left(u_*\right).
\end{align*}

\bigskip

As above, the threshold $c_\varepsilon$ is given implicitly as the solution of
the equation~\eqref{thrC}. In the following table we give some values of
$c_\varepsilon$ obtained using numerical simulations.

\begin{center}
\begin{tabular}{| l | c | c | c | c | c | c |}
\hline
$\varepsilon$&0.01&0.02&0.03&0.04&0.05&0.1
\\\hline
$c_\varepsilon$&13.692&11.224&9.803&8.806&8.042&5.719
\\\hline
\end{tabular}
\end{center}

These thresholds are obtained by simulating $M=10^7$ trajectories on $[0,1]$
of a standard Wiener process, calculating for each of them the quantity
$\Gamma(W)$ and taking $(1-\varepsilon)M$-th greatest between them.

\section{Comparison of the Tests}

  Remind that all these three tests $\phi_T^*,\bar\phi_T $ and $\hat\phi_T$ in
regular (LAN) case are asymptotically equivalent to the Neyman-Pearson test
$\phi_{u,T}^\circ $ (with known alternative $u$) and hence are asymptotically
uniformly most powerful. In our singular situation all of them have different
asymptotic behavior and therefore it is interesting to compare their limit
power functions
\begin{align*}
\beta^*\left(u\right)&=\Pb_u\left\{\Delta \left(Y_u\right)>a_\varepsilon
\right\},\qquad \bar \beta\left(u\right)=\Pb_u\left\{\frac{\Delta
\left(Y_u\right)}{\sqrt{2{\rm J}\left(Y_u\right)}}>b_\varepsilon
\right\},\\
\hat\beta\left(u\right)&=\Pb_u\left\{\frac{\Delta
\left(Y_u\right)}{{\rm J}\left(Y_u\right)}>c_\varepsilon \right\},\qquad
\beta^\circ\left(u\right)=\Pb_u\left\{u\Delta \left(Y_u\right)-\frac{u^2}{2}
{\rm J}\left(Y_u\right)>d_\varepsilon \right\}
\end{align*}
of course, under condition that all of them belong to $ {\cal K}_\varepsilon
$. Our goal is to compare these quantities for the large values of $u$.

We have to study the distribution of the vector $\left(\Delta
\left(Y_u\right),{\rm J}\left(Y_u\right) \right)$, where
$$
\Delta \left(Y_u\right)=-\int_{0}^{1}Y_u\left(s\right)\;{\rm
d}Y_u\left(s\right),\qquad {\rm
J}\left(Y_u\right)=\int_{0}^{1}Y_u\left(s\right)^2\;{\rm d}s,
$$
where $Y_u $ is solution of the equation
$$
{\rm d}Y_u\left(s\right)=-u\;Y_u\left(s\right)\;{\rm d}s+{\rm
d}W\left(s\right),\qquad Y_u\left(0\right)=0 ,\qquad 0\leq s\leq 1.
$$

Let us introduce the stochastic process $y_v=\sqrt{u}\;Y_u\left(
\frac{v}{u}\right),0\leq v\leq u$ (this transformation was introduced by
Luschgy \cite{Lus2}). Then we can write
$$
{\rm d}y_v=-y_v\;{\rm d}v+{\rm d}w_v,\qquad y_0=0,\qquad 0\leq v\leq u,
$$
where $w_v =\sqrt{u}\;W\left(\frac{v}{u}\right)$ is a Wiener process and
$$
\Delta \left(Y_u\right)=-u^{-1}\int_{0}^{u}y_v\;{\rm d}y_v\equiv\frac{\Delta_u
}{u} ,
\qquad {\rm J}\left(Y_u\right)= u^{-2}\int_{0}^{u}y_v^2\;{\rm
d}v\equiv\frac{{\rm J}_u }{u^2}.
$$
in obvious notation. Further, the process $y_v$ is ergodic with the density of
the invariant law $f\left(y\right)= e^{-y^2}/\sqrt{\pi} $. Hence ${\rm
J}_u\rightarrow \infty $ and
$$
\frac{1}{u}\int_{0}^{u}y_v^2\;{\rm d}v\longrightarrow \frac{1}{2}.
$$
Note that the distribution of the process $y_v$ does not depend on $u$.

The constant $d_\varepsilon =d_\varepsilon \left(u\right)$ because it is
defined by the equation
$$
\Pb_0\left\{u\Delta \left(W\right)-\frac{u^2}{2}
{\rm J}\left(W\right)>d_\varepsilon \right\} =\varepsilon .
$$
For the large values of $u$ this constant can be approximated as follows. We
have (under hypothesis ${\cal H}_0$) as $u\rightarrow \infty $
\begin{align*}
 &\Pb_0\left\{u\Delta \left(W\right)-\frac{u^2}{2}{\rm
J}\left(W\right)>d_\varepsilon \left(u\right)  \right\}=\\
&\qquad =\Pb_0\left\{\int_{0}^{1}W\left(s\right)^2{\rm d}s<-\frac{2d_\varepsilon
\left(u\right)}{u^2} + \frac{2\Delta\left(W\right)}{u}\right\} \longrightarrow\\
&\qquad \longrightarrow
\Pb_0\left\{\int_{0}^{1}W\left(s\right)^2{\rm d}s<e_\varepsilon \right\}
=\varepsilon,
\end{align*}
 where the constant $e_\varepsilon $ is defined by the last
equality. For example, if we take $\varepsilon =0,05$ then the numerical
simulation gives us the value $e_{0,05}=0,056$. Therefore
$d_\varepsilon\left(u\right)=-0,5\,e_\varepsilon
\,u^2\left(1+o\left(1\right)\right) $. If we suppose that $\varepsilon $ is
small and try to solve the equation
$$
\int_{0}^{e_\varepsilon }f_{{\rm J} }\left(x\right)\,{\rm d}x=\varepsilon
$$
where $f_{{\rm J} }\left(x\right) $ is the density function of the integral ${\rm
J}\left(W\right) $, then we can easily see that $f_{{\rm J} }\left(0\right)=0 $
and all its derivatives $f_{{\rm J} }^{\left(k\right)}\left(0\right)=0, k=1,2,
\ldots$. Hence to see an approximative solution we need to calculate the large deviation
probability of the following form (below $r=s/\sqrt{e_\varepsilon },
E=e_\varepsilon ^{-1/2}\rightarrow \infty $).
\begin{align*}
\Pb_0\left\{e_\varepsilon ^{-1} \int_{0}^{1}W\left(s\right)^2{\rm
d}s<1\right\}= \Pb_0\left\{\int_{0}^{E}W\left(r\right)^2{\rm
d}r<1\right\} .
\end{align*}
Below we put $d_\varepsilon\left(u\right)=-0,5\,e_\varepsilon \,u^2 $.

We have the relations
\begin{align*}
\beta^*\left(u\right)&=\Pb\left\{\Delta_u>u\,a_\varepsilon
\right\}=\Pb\left\{\int_{0}^{u}y_v\,{\rm d}w_v<{\rm J}_u-a_\varepsilon\,u\right\},\\
 \bar \beta\left(u\right)&=\Pb\left\{\frac{\Delta_u
}{\sqrt{2{\rm J}_u}}>b_\varepsilon
\right\}=\Pb\left\{\int_{0}^{u}y_v\,{\rm d}w_v<{\rm
J}_u-b_\varepsilon\,\sqrt{2{\rm J}_u }\right\},\\
\hat\beta\left(u\right)&=\Pb\left\{\frac{\Delta_u
}{{\rm J}_u}>\frac{c_\varepsilon }{u}\right\}=\Pb\left\{\int_{0}^{u}y_v\,{\rm
d}w_v<{\rm
J}_u-\frac{c_\varepsilon}{u}\,{\rm J}_u \right\},\\
\beta^\circ\left(u\right)&=\Pb\left\{\Delta_u-\frac{{\rm J}_u}{2}
>d_\varepsilon \right\}=\Pb\left\{\int_{0}^{u}y_v\,{\rm d}w_v<\frac{1}{2}{\rm
J}_u+\frac{e_\varepsilon}{2}\,u^2 \right\} .
\end{align*}

Therefore the large values of $u$ (${\rm J}_u\sim u/2 $)
\begin{align*}
&\frac{1}{2}{\rm J}_u+\frac{e_\varepsilon}{2}\,u^2 > {\rm
J}_u-\frac{c_\varepsilon}{u}\,{\rm J}_u>{\rm
J}_u-b_\varepsilon\,\sqrt{2{\rm J}_u }>{\rm J}_u-a_\varepsilon\,u,
\end{align*}
and finally
$$
\beta^*\left(u\right)<\bar
\beta\left(u\right)<\hat\beta\left(u\right)<\beta^\circ\left(u\right).
$$
These inequalities are in accord with \cite{Swen97}.

Note that for small values of $\varepsilon$ the constant $a_\varepsilon $ is
close to 0,5 (e.g. $a_{0,05}=0,498$, $a_{0,01}=0,49992$) and in this
asymptotics the power of score-function test is
$$
\beta^*\left(u\right) =\Pb\left\{\int_{0}^{u}y_v\,{\rm d}w_v<
\left(0,5-a_{\varepsilon }\right)\,u\left(1+o\left(1\right)\right)\right\}.
$$
Hence one can expect that in this case the score-function test has essentially
smaller power than the others.

\bigskip

Now let us turn to numerical simulations of the limiting power functions. We
aim to obtain the limiting power functions of all the three tests, as well as
the Neyman-Pearson envelope, for the moderate values of $u$ ($u\leq15$).

Note that for the score function test $\beta^*(u)$ can be computed
directly using~\eqref{10}.  However the limiting power functions of the
likelihood ratio and of the Wald tests are written as probabilities of some
events related to Ornstein-Uhlenbeck process and can be obtained using
numerical simulations.

For the likelihood ratio test we have
$$
\bar\beta\left(u\right)=
\Ex_u \1_{\left\{\Lambda(Y_u)>b_\varepsilon\right\}}=
\Ex_0 Z\left(u\right)\1_{\left\{\Lambda(W)>b_\varepsilon\right\}}
$$
where
$$
Z\left(u\right)=\exp\left\{u\Delta(W)-\frac{u^2}{2}\,{\rm J}(W)\right\}.
$$

So we simulate $M=10^7$ trajectories $W_j=\left\{W_j(s),\ 0\leq s\leq
1\right\}$, $j=1,\ldots,M$ of a standard Wiener process and calculate for each
of them the quantities $\Delta_j=\Delta(W_j)$, ${\rm J}_j={\rm J}(W_j)$,
$\Lambda_j=\Delta_j/{\rm J}_j$ and (for each value of $u$)
$Z_j\left(u\right)=\exp\left\{u\Delta_j-\frac{u^2}{2}\,{\rm
J}_j\right\}$. Then we calculate the empirical mean
$$
\frac1M\sum_{j=1}^{M} Z_j\left(u\right)\1_{\left\{\Lambda_j>b_\varepsilon\right\}}
\approx\bar\beta\left(u\right).
$$

For the Wald test we have similarly
$$
\frac1M\sum_{j=1}^{M} Z_j\left(u\right)\1_{\left\{\Gamma_j>c_\varepsilon\right\}}
\approx\hat\beta\left(u\right)
$$
where $\Gamma_j=\Delta_j/\sqrt{2\,{\rm J}_j}$.

Finally, in order to compute the Neyman-Pearson envelope, we first approximate
(for each value of $u$) the quantity $d_\varepsilon=d_\varepsilon(u)$ by the
$(1-\varepsilon)M$-th greatest between the quantities $\ln Z_j(u)$, and then
calculate
$$
\frac1M\sum_{j=1}^{M} Z_j\left(u\right)\1_{\left\{\ln
Z_j(u)>d_\varepsilon(u)\right\}}\approx\beta^\circ\left(u\right).
$$

The results of these simulations for $\varepsilon=0.05$ are presented in
Figure~2.

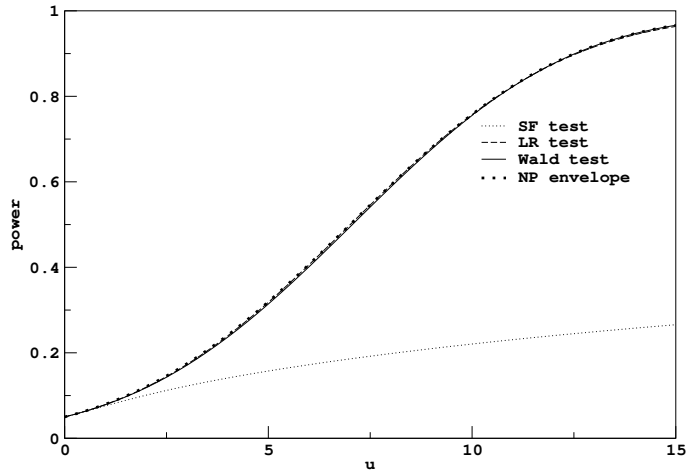
\onepic{SCpow1.eps}{Fig.~2: Limiting powers for $\varepsilon=0.05$}

Let us note here that in this case the power functions of the likelihood ratio
test and of the Wald test are indistinguishable (from the point of view of
numerical simulations) from the Neyman-Pearson envelope. This quite surprising
fact was already mentioned by Eliott {\sl at al.}  \cite{Eli}, who showed the
similar pictures having $2\cdot 10^3$ simulations. As we see from Figure 2,
with $10^7$ simulations the curves are still indistinguishable. The situation
is however different for bigger values of $\varepsilon$. The results of
simulations for $\varepsilon=0.01$, $0.05$, $0.25$ and $0.5$ are presented in
Figure~3.

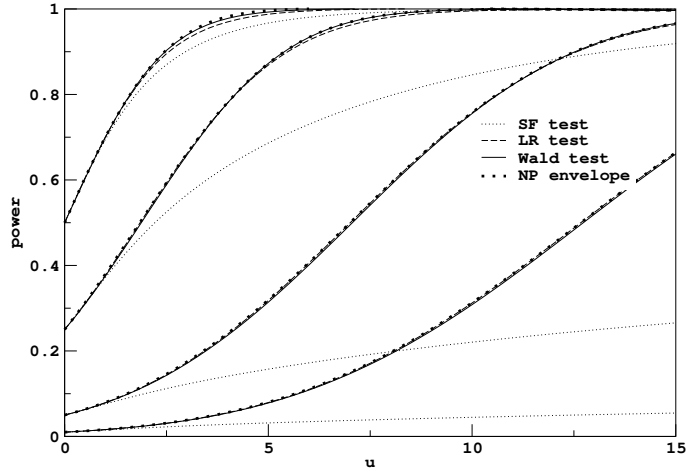
\onepic{SCpow3.eps}{Fig.~3: Limiting powers for different values of
$\varepsilon $}

One can note that for big values of $\varepsilon$ (e.g. $\varepsilon=0.5$) the
powers became more distinguishable, and that the asymptotically established
ordering of the tests holds already for these moderate values of $u$. Note
also that for the small values of $\varepsilon$ (e.g. $\varepsilon=0.01$ and
$0.05$) the curve of score-function test is essentially lower as expected.

\section{Discussion}

\noindent {\bf Remark 1.} Note that alternatives $u=u_T\rightarrow \infty $
with $\vartheta _{u_T}\rightarrow 0$ are local but not contigous. That means
that the corresponding sequences of measures
$\left(\Pb^{\left(T\right)}_{\vartheta
_{u_T}},\Pb^{\left(T\right)}_{0}\right), T\rightarrow \infty $ are not
contigous. Particularly, the second integral in the likelihood ratio formula
tends to infinity:
$$
\int_{0}^{T}\left[\psi \left(\vartheta _{u_T}
\left(S_*t-X_{t-}\right)\right)-1- \ln \psi \left(\vartheta _{u_T}
\left(S_*t-X_{t-}\right)\right)\right]\;S_*\,{\rm d}t\longrightarrow \infty .
$$
In such situation the power function of any reasonable test tends to 1 and to
compare tests we have to use, say, the large deviation principle.  For
example, the likelihood ratio test $\phi_\zT^* $ is consistent for the {\sl
local far alternatives} $\vartheta =\frac{v}{\sqrt{S_*T}}, v\in \left[\nu
,V\right]$ where $0<\nu<V< \infty $. Indeed, under mild regularity conditions
we can write
\begin{align*}
&\Ex_v\phi_\zT^*\left(X^T\right)= \Pb_0\left\{\sup_{\nu <v<V} L\left(\frac{v}{
\sqrt{S_*T}},X^T\right)>c_\varepsilon \right\}=\\
&\qquad = \Pb_0\left\{\sup_{\nu
<v<V}\left[\sqrt{S_*T}\int_{0}^{1}\ln\psi\left(vW_T\left(s\right)\right){\rm
d}W_T\left(s\right)-\right.\right.\\
&\qquad \quad
\left.\left.
-S_*T\int_{0}^{1}\left[\psi\left(vW_T\left(s\right)\right)-1-
\ln\psi\left(vW_T\left(s\right)\right)\right]
{\rm d}s\right]>\ln c_\varepsilon \right\}=\\
&=\Pb_0\left\{\sup_{\nu
<v<V}\left[\frac{1}{\sqrt{S_*T}}\int_{0}^{1}\ln\psi\left(vW_T\left(s\right)\right){\rm
d}W_T\left(s\right)-\right.\right.\\
&\qquad \quad
\left.\left.
-\int_{0}^{1}\left[\psi\left(vW_T\left(s\right)\right)-1-
\ln\psi\left(vW_T\left(s\right)\right)\right]
{\rm d}s\right]>\frac{\ln c_\varepsilon}{S_*T} \right\}\longrightarrow \\
&\longrightarrow \Pb\left\{\inf_{\nu
<v<V} \int_{0}^{1}\left[\psi\left(vW\left(s\right)\right)-1-
\ln\psi\left(vW\left(s\right)\right)\right]{\rm d}s>0\right\}=1
\end{align*}
because the function $g\left(y\right)=y-1-\ln y> 0$ for $y\neq 1$ and
$g\left(y\right)=0$ iff $y=1$.

\bigskip

\noindent
{\bf Remark 2.} Note, that we can construct asymptotically uniformly most
powerful test if we change the statement of the problem in the following
way. Let us fix some $D>0$ and introduce the stopping time
$$
\tau _D=\inf \left\{\tau :\;\int_{0}^{\tau } \left(S_*t-X_t\right)^2\;S_* \;{\rm
d}t\geq  D^2\right\}.
$$

Then  we consider the problem of testing hypotheses
\begin{align*}
{\scr H}_0\quad :\quad\qquad  S\left(t,X_t \right)&=S_*,\\
{\scr H}_1\quad :\quad\qquad S\left(t,X_t \right)&=S_*\; \psi\left(\vartheta
_D\left[S_*t-X_t\right]\right), \quad \vartheta _D=\frac{u}{\dot\psi\left(0\right)D}>0
\end{align*}
by observations $X^{\tau _D}=\left\{X_t,0\leq t\leq \tau _D\right\}$ in the
asymptotics  $D\rightarrow \infty $.  Now the likelihood ratio
$Z_{\tau_{D}}\left(u\right)=L\left(\frac{u}{\dot\psi\left(0\right)D},X^D\right)$
will be LAN:
$$
Z_{\tau_{D}}\left(u\right)\Longrightarrow \exp\left\{u\;\zeta
-\frac{u^2}{2} \right\},\qquad \zeta \sim
{\cal N}\left(0,1\right)
$$
and the test $\hat\phi_{\tau_{D}}=\1\zs{\left\{\Delta
_{\tau_{D}}\left(X^{\tau_{D}}\right)>z_\varepsilon \right\}}$ where
$$
\Delta_{\tau_{D}}\left(X^{\tau_{D}}\right)=\frac{1}{D}\int_{0}^{\tau
_D}\left(S_*t-X_{t-}\right)\;\left[{\rm d}X_t-S_*{\rm d}t\right]
$$
is locally asymptotically uniformly most powerful.

The proof follows from the central limit theorem for stochastic integrals and
the standard arguments (for  LAN families).

\bigskip
\noindent {\bf Remark 3.} Note that these problems of hypotheses testing are
similar to the corresponding problems of hypotheses testing for diffusion
processes. In particular, let the observed process $X^T=\left\{X_t,0\leq t\leq
T \right\}$ be diffusion
$$
{\rm d}X_t= \psi \left(-\vartheta_T\,X_t\right)\;{\rm d}t+\sigma \,{\rm d}W_t,\quad
X_0=0, \quad 0\leq t\leq T ,
$$
where the function $\psi\left(0 \right)$=0, is continuously differentiable at
the point $0$ and $\dot{\psi}\left(0 \right)>0$.  If we consider two
hypotheses: $\vartheta =0 $ and $\vartheta >0 $ then the reparametrization
$$
\vartheta _T=\frac{u\,\sigma}{\dot{\psi}\left(0 \right)\;T }
$$
provides local contigous  alternatives, i.e., the log-likelihood ratio in the problem
\begin{align*}
{\scr H}_0&:\quad\qquad \quad\qquad u=0,\\
{\scr H}_1&:\quad\qquad \quad\qquad u>0.
\end{align*}
 has the  limit:
$$
\ln L\left(\frac{u\,\sigma}{\dot{\psi}\left(0 \right)\;T
},X^T\right)\Longrightarrow -u\int_{0}^{1}W\left(s\right) {\rm
d}W\left(s\right)-\frac{u^2}{2}\int_{0}^{1}W\left(s\right)^2  {\rm d}s.
$$
The score function test based on the statistic
$$
\Delta _T^*\left(X^T\right)= -\frac{1}{T}\int_{0}^{T}X_{t}\;{\rm d}X_t,
$$
 the likelihood ratio test and the Wald test have the same asymptotic properties as those
described in Theorems 1, 2  and 3 above.

For example, if $\psi\left(x\right)=x$, then we have the  Wiener process (under
hypothesis ${\scr H}_0$) against ergodic  Ornstein-Uhlenbeck
process under alternative ${\scr H}_1$.

\bigskip
\noindent {\bf Remark 4.} We supposed above that the derivative of the
function $\psi\left(x \right)$ at the point $x=0$ is not equal to 0, but
sometimes it can be interesting to study the score function and the likelihood
ratio test in the situations when the first $k-1$ derivatives with $k\geq 2$ are
null.

Let us consider a self-correcting process $X^T=\left\{X_t,0\leq t\leq
T\right\}$ with intensity function $S_*\psi\left(\vartheta
\left(S_*t-X_t\right) \right)$ such that $\psi\left(0 \right)=1$
$\dot\psi\left(0 \right)=0$ and $\ddot\psi\left(\cdot \right)\neq 0$
$\left(k=2\right)$. In this case the modifications have to be the
following. Suppose that $\ddot\psi\left(0 \right)> 0$. To have LAQ family at
the point $\vartheta =0$ we chose the reparametrization $\vartheta =\vartheta
_u$
$$
\vartheta_u
=\sqrt{\frac{2\,u}{\ddot\psi\left(0\right)}}\;\left(S_*T\right)^{-3/4},
$$
which provides the limit
$$
\ln L\left(\vartheta_u,X^T\right)\Longrightarrow
u\;\int_{0}^{1}W\left(s\right)^2\;{\rm
d}W\left(s\right)-\frac{u^2}{2}\;\int_{0}^{1}W\left(s\right)^4\;{\rm d}s
$$

Then in the hypotheses testing problem
\begin{align*}
{\scr H}_0&:\quad\qquad \quad\qquad u=0,\\
{\scr H}_1&:\quad\qquad \quad\qquad u>0
\end{align*}
 the score function test $\hat\psi\left(X^T\right)=\1_{\left\{\Delta
 _\zT\left(X^T\right)>c_\varepsilon  \right\}}$ is based on the statistic
$$
\Delta _\zT\left(X^T\right)=\frac{1}{\left(S_*T\right)^{3/2}}\int_{0}^{T}
\left(S_*t-X_t\right)^2\;\left[{\rm d}X_t-S_*{\rm d}t\right].
$$
It is easy to see that under ${\scr H}_0$
$$
\Delta _\zT\left(X^T\right)\Longrightarrow
\frac{W\left(1\right)^3}{3}-\int_{0}^{1}W\left(s\right)\;{\rm d}s .
$$
Hence to chose the threshold $c_\varepsilon $ we have to solve the following equation
$$
\frac{4}{3}\iint_{x^3-y>3c}\exp\left\{-2x^2+2xy-\frac{2}{3}y^2\right\}\;\;{\rm
d}x\,{\rm d}y=\varepsilon
$$
because $\left(W\left(1\right),3\int_{0}^{1}W\left(s\right)\;{\rm d}s\right) $
is Gaussian vector.

The cases $k>2$ can be treated in a similar way.

\end{document}